\documentclass{article}
\usepackage{graphicx}
\usepackage{amsmath}
\usepackage{amssymb}
\usepackage{hhline}
\usepackage{a4wide}
\usepackage{mathrsfs}
\newtheorem{theorem}{Theorem}

\newtheorem{theoremc}{Theorem}

\newtheorem{prop}[theorem]{Proposition}
\newtheorem{rk}[theoremc]{Remark}
\newenvironment{proof}[1][Proof]{\textbf{#1.} }{\qed \vspace{5pt}}

\renewcommand\l{\lambda}

\newcommand\vp{\varphi}

\newcommand\bib[1]{\bibitem[#1]{#1}}

\newcommand\abz{\hspace{13.5pt}}
\newcommand\qed{\phantom{\underline{y}}\hfill\hfill$\square$\vspace{3pt}}

\newcommand\C{{\mathbb C}}
\newcommand\Cc{{\let\mathcal\mathscr\mathcal C}}
\newcommand\D{{\mathcal D}}
\newcommand\op[1]{\mathop{\rm #1}\nolimits}
\newcommand\E{{\mathcal E}}

\newcommand\R{{\mathbb R}}
\newcommand\hps{\hskip-16pt . \hskip2pt}
\newcommand\hpss{\hskip-13.5pt . \hskip2pt}

\newcommand\z{\sigma}
\newcommand\p{\partial}

\newcommand\1{{\bf 1}}

\DeclareFontFamily{U}{wncyr}{}
\DeclareFontShape{U}{wncyr}{m}{it}{%
  <5><6><7><8><9>gen*wncyi%
  <10><10.95><12><14.4><17.28><20.74><24.88>wncyi10}{}
\DeclareSymbolFont{MathRussLetters}{U}{wncyr}{m}{it}
\DeclareMathSymbol{\re}{\mathalpha}{MathRussLetters}{3}

\newcommand{\weg}[1]{}

\makeatletter
\renewcommand{\@oddhead}{\hfil Reductions via Compatibility\hfil}
\renewcommand{\@evenhead}{\hfil Boris Kruglikov\hfil}
\renewcommand{\@begintheorem}[2]{\begin{trivlist}\it
 \item[\hspace{\labelsep}{\bf #1\ #2.}]}
\renewcommand{\@endtheorem}{\end{trivlist}}
\makeatother

\begin{document}

 \title{Symmetry approaches for reductions of PDEs,\\
 differential constraints and Lagrange-Charpit method}
 \author{Boris Kruglikov}
 \date{}
 \maketitle

 \vspace{-14.5pt}
 \begin{abstract}
Many methods for reducing and simplifying differential equations are known. They provide various generalizations of the original symmetry approach of Sophus Lie. Plenty of relations between them have been noticed and in this note a unifying approach will be discussed. It is rather close to the differential constraint method, but we make this rigorous basing on recent advances in compatibility theory of non-linear overdetermined systems and homological methods for PDEs.
 \footnote{MSC numbers: 34A05, 35N10; 58A20, 35A30.\\
 Keywords: Compatibility, differential constraint, symmetry, reduction, multi-bracket, solvability.}%
 \end{abstract}

\section*{Introduction to reduction methods}\label{S0}

 \hspace{13.5pt}
Exact solutions play an important role in investigation of equations of mathematical physics. Even though usually\footnote{Except for the important class of completely integrable equations.}
it's only a small portion of solutions one can find explicitly, they represent significant physical phenomena and are often applied in practice. The first explicit methods have origins in the works of Euler, Lagrange, Cauchy, Hamilton, Jacobi and Monge. But the most essential contribution was made by Sophus Lie, whose continuous transformation groups changed the landscape of the theory of differential equations. These transformations are called nowadays symmetries\footnote{They are distinguished to be point or contact symmetries, though the latter received less circulation.} and can be computed with the help of most popular symbolic software packages.

If the symmetry algebra is known, is sufficiently large and has an appropriate Lie structure, then the given differential equation (or system) can be integrated completely. This works both for ODEs and PDEs. But usually the algebra is not big\footnote{For linear system it is  huge, but the bulk consists of shifts by solutions and the quotient is usually moderately small.}, so other methods should apply.

There is a number of classical methods, like intermediate integrals, differential and nonlocal substitutions (e.g. Laplace, Backlund and Darboux transformations) and conservation laws. In more recent time a lot of other methods for obtaining exact solutions were invented. The whole list is large, but we name a few: higher and inner (=internal) symmetries \cite{KLV}, recursion operators, nonlocal and ghost symmetries, nonclassical symmetries and direct reduction, inverse scattering, bi-linear Hirotha method, Sato theory and so on, see e.g. \cite{BK,CRC,SCL}.

Many interrelations have been made explicit, but the unifying theory does not exist so far. In this paper I will try to present a moderately general viewpoint for all local theories (nonlocal ones follow the same scheme, but additionally need a covering, i.e. another system of DE implying the given one as a differential corollary). I will also try to explain why it is cumbersome to encompass all features into the framework of one theory leaving it reasonably efficient.

 \weg{
There are two poles in investigation of PDEs: one is based on the functional analysis approach, which proves existence and even gives methods to further numerical establishing and investigations of general solutions of linear and non-linear differential equations. Usually in such cases quite arbitrary functions enter the set-up as initial or boundary values and there is no hope to find an exact form of solution except some particular cases, where a nice Fourier or integral representation exists. However in many such systems the results are sensitive to functional spaces one chooses, and they often differ from $C^\infty$ spaces (except for some elliptic kind linear problems, where we can finally settle to smooth functions), thus physical phenomena like blow-ups, shocks etc being dependent on the sort of completion.

Another pole is the mathematical physics approach, where we wish to get all solutions in more or less explicit way, i.e. involving certain known functions and operators of inversion, quadrature and maybe bit more general. This however is rarely possible, though the equations allowing such treatment plays a crucial role. An intermediate step consists of equations, which are not completely integrable, but for which a large set of explicit solutions can be obtained...
 }

\subsection{\hpss Differential constraints.}\label{S01}

 \abz
This method is pretty old (it dates back at least to the time of Lagrange) and is quite obvious: given a system $\E=\{F_1=\dots=F_r=0\}$ one overdetermines it with other differential equations $\mathcal{G}=\{G_1=\dots=G_s=0\}$ and investigates the joint system $\E\cap\mathcal{G}=\{F_i=0,G_j=0\}$. Clearly $\op{Sol}(\E\cap\mathcal{G})\subset\op{Sol}(\E)$ (here and in what follows by $\op{Sol}(\E)$ we understand local or even formal solutions, though global will work as well).

 \begin{rk}\label{rk1}
This procedure is the direct analog of the standard tool in algebraic geometry: In order to investigate an algebraic variety $X\subset\mathbb{P}^N$ it is bi-sected by another variety $Y$ and the totality of intersections $X\cap Y$ uncover $X$ (most often $Y$ is taken in the simplest way, like an affine plane $\mathbb{P}^k\subset\mathbb{P}^N$, but sometimes more complicated variety is needed). In the case of differential equations the section $\op{Sol}(\mathcal{G})$ may be even more complicated than the original $\op{Sol}(\E)$, but the task is to make $\op{Sol}(\E\cap\mathcal{G})$ tractable.
 \end{rk}

We will suppose that the original system $\mathcal{E}$ is compatible (=formally integrable). Let us call the constraint $\mathcal{G}$ compatible differential constraint if the joint system $\E\cap\mathcal{G}$ is compatible as well (NB: $\mathcal{G}$ is not required to be compatible!).

This method was introduced into practice by Yanenko \cite{Y} and also
as a generalization of S.Lie approach by Ovsyannikov \cite{Ov}, see also Bluman and Cole \cite{BC}, though relations between these were unexplored for some time \cite{Me}. Re-invented by the name of side conditions in \cite{OR}, the method was much popularized later by Olver \cite{Ol}. In fact, it was in this latter papers that the most cited statement was made (we change it a bit and state the claim in the ultimate form):

 \begin{quote}
\hspace{-12pt}{\it All reduction methods are partial cases of the method of compatible differential constraints.}\hspace{-10pt}
 \end{quote}

In particular it was shown in \cite{Ol} that many known methods such as classical and non-classical symmetry reductions, direct method of Clarkson-Kruskal, Galaktionov's separation method and some others are covered by the compatible constraints.

This is an important statement and some remarks should be of order.

{\bf 1.} The result is demonstrated for a scalar single second order PDE on the plane only and even though it was claimed that "the results can be straightforwardly generalized to arbitrary systems", it is neither clear to what degree (see 2 below which shows one must be precise) nor how (direct generalization of this proof will lead to the trouble even with a single scalar PDE of order 3 on $\R^3$; consult \cite{KL$_1$,KL$_4$} for a similar but rigorous statements to appreciate which technique should be used for quite general, though still not arbitrary systems).

{\bf 2.} The statements in \cite{OR,Ol} are rather vague, but without specific genericity conditions the results turn wrong. For instance the following is a combination of Proposition 1 and Theorem 2 from \cite{Ol}: {\it A first order differential constraint $G=0$ reduces the second order PDE $F=0$ to an ODE iff the system $F=0,G=0$ is compatible\/}. The trivial counter-example to this is $F=H\circ G$, which is compatible with $G$ for any differential operator $H$ (if the operators are non-linear we should rather use linearization in the composition: $F=\ell_H\circ G$).

{\bf 3.} The situation becomes worse in higher order. For instance Darboux semi-integrable systems can be treated as compatible pairs $F=0,G=0$ of scalar PDEs on the plane with 1 common characteristic, which surely declines possibility of direct reduction to an ODE. For higher order systems on $\R^n$ the situation becomes completely untractable if no restrictions are imposed (the whole algebraic singularity theory imbeds into this problem!).

Lacking a general proof and due to these remarks we regard the above claim as a conjecture, whose precise form should be specified. The proof will be given in the subsequent sections.

\subsection{\hpss General statement.}\label{S02}

 \abz
Reductions of systems of PDEs to systems of ODEs lead to final dimensional families of solutions. Actually this can be taken as a definition of reduction: Any stratified manifold $\mathcal{S}\subset\op{Sol}(\E)$ will be called a reduction to ODEs. This in particular covers the case when ODEs have singular solutions because these should be considered as components of $\mathcal{S}$ of different dimensions. The claim formalizing the one of the previous section is now the following:

 \begin{theorem}\label{thm1}
Any reduction $\mathcal{S}\subset\op{Sol}(\E)$ can be obtained via the method of compatible differential constraint.
 \end{theorem}

This statement is not too rigorous either. The main problem is the space $\op{Sol}(\E)$ that is not defined due to the lack of existence and uniqueness theorems in PDE theory. Some solutions to this exist. One is to consider analytic systems and analytic solutions or to boil everything down to the formal theory (this requires a lot of formalism but is possible). Another one is to restrict to classes of smooth PDEs, which have existence and uniqueness theorems (together with regularity theorems so that a-posteriori we do not need to leave the class of smooth functions). We would prefer however to consider the stratified manifold $\mathcal{S}$ abstractly and claim that all its points are (loc.) solutions of $\E$. There are standard ways to impose the topology on it.

\medskip

 \begin{proof}
We can fix a setup so that $\mathcal{S}\subset C^\infty(\pi)$ for a bundle $\pi:E\to M$. It is possible to do globally due to Whitney theorem for stratified manifolds, but local representation will suffice for us because we'll be patching them later on.

Consider all (non-linear) differential operators $G\in\op{diff}(\pi,\1)$ such that $G[u]=0$ for any $u\in\mathcal{S}$. This is a set $Y(\mathcal{S})\supset Y(\E)=\langle F_1,\dots F_r\rangle$. The space $Y(\mathcal{S})$ is a differentially closed algebra. Indeed if $G_i\in Y(\mathcal{S})$ vanish on section $u\in\mathcal{S}$, so do their linear combinations and left-compositions with other differential operators.

Now we claim that $Y(\mathcal{S})$ is finitely generated over $\op{diff}(\1,\1)$. In formal theory this is so due to Noetherian property (one can consult \cite{KL$_5$} for the framework of formal theory for non-linear operators, when modules should be considered over the algebra of $\Cc$-differential operators).

In analytic context this is generally wrong, but we can assume that $Y(\mathcal{S})$ is a prime ideal and then this becomes the fundamental theorem of Ritt-Kolchin theory \cite{Ri,Ko} (we refer to \cite{B,Ma} for the theory of D-modules and applications to our context).

In smooth case one should act as follows: a neighborhood of a stratum in $\mathcal{S}$ is parametrized by constants $c_1,\dots,c_d$. Taking the space of all parametrized derivatives $D_\z u(c_1,\dots,c_d)$ the implicit function theorem makes it possible to specify the defining functions of the $d$-dimensional submanifold with parameters $c_i$. This finite set can always be specified to begin with $F_i$, so it is $F_1,\dots,F_r,G_1,\dots,G_s$. The second part of this collection defines $\mathcal{G}$.

Note however that in any of the above approaches it is not possible to distinguish the whole manifold $\mathcal{S}$ at once. Overdeterminations $\mathcal{G}$ play the role of local charts and as in differential geometry we need a finite or countable collection of such. This is especially clear with application of implicit function theorem, but is necessary for other cases as well. Finally we patch the local charts and so $\mathcal{S}$ is defined by the compatible differential constraint $\mathcal{G}$.
 \end{proof}

Even though on the level of ideas this proof is clear, the reduction methods have quite precise implementation techniques and algorithms. A prompt objection to the above is that it is not very constructive. But this will be compensated in the next section, where we explain how standard methods imbeds into compatibility theory and make specifications.

 \begin{rk}\label{rk2}
The above reduction to ODE, in light of the analogy from Remark \ref{rk1}, is similar to the section of $X$ by transversal varieties $Y$ of complimentary dimension, $\op{codim}Y=\op{dim}X$, when $\dim(X\cap Y)=0$. Note however that with the above approach we always obtain not only solutions $\mathcal{S}$, but more generally its D-Zarissky closure $\bar{\mathcal{S}}$.

Higher-dimensional strata of $\op{Sol}(\E)$ are parametrized by $\z$ functions of $p$ variables ($\z$ is the functional rank and $p$ the functional dimension of $\op{Sol}(\E)$, we refer to \cite{KL$_3$} for the discussion). Thus we can also make a more general reduction of $p$ or keep it and reduce $\z$.

Note however that several small reductions are not equivalent to one big, because several steps can be impossible, see \cite{Cl}\S6.5. This is known already for the classical symmetry reduction of ODEs, which do not respect the reductive algebra structure of the symmetry algebra \cite{Ly}.
 \end{rk}

\subsection{\hpss Discussion.}\label{S03}

 \abz
Let us make distinction between compatible differential constraints, solvable differential constraints and just differential constraints. The former are based on some explicit compatibility criteria which, provided the type of symbolic syzygy (or type of singularity for characteristic variety) is fixed\footnote{In most geometric approaches to differential equations this is a crucial requirement of regularity, so that the type (symbols, characteristics, Spencer $\delta$-cohomology etc) does not change along the equation considered as a submanifold in jets. Equivalently one can require that the set of singular points is meager and can be removed. For analytic D-modules it can be proved that the set of these point has Zarisky-open complement and thus is negligible.}, can be written explicitly, see the next section.

On the contrary there do not exist general solvability criteria for overdetermined systems. All of them are based on the prolongation-projection method, which consists in adding compatibility conditions and removing irregular points. Thus it is neither explicit nor simple. Indeed the computational cost of this process can be quite large\footnote{Even polynomial Gr\oe bner basis has double exponential complexity, its differential analogs (\cite{Hu} and references therein) exhibit no decrease, see e.g. \cite{Kr} for this question in relation with generalized Lagrange-Charpit method.}. Thus there are no big difference between solvable differential constraints and just differential constraints: One starts with an arbitrary overdetermination and checks consistency getting at the end formally integrable system or contradiction (Cartan's formulation for inconsistent equation \cite{Ca}) as Cartan-Kuranishi theorem states \cite{Ra}.

In this latter form the method of differential constraints is equivalent to the method of ansatzes. It is a matter of taste what to choose then. When we join an additional PDE, it is usually easier to check compatibility conditions. On the other hand, an ansatz is ready to provide the form of solution, while with differential constraint one should first solve the auxiliary equations.

\section{\hps Details of the method}\label{S1}

 \abz
Since compatible differential constraints play the central role in reduction of differential equations it is important to have efficient compatibility criteria. The method of E.Cartan, Jane and Riquier as well as Spencer theory (see \cite{KLV,KL$_5$} and references therein) provide theoretical base and algorithm for establishing compatibility conditions, but they give no formulas. The latter were established for a wide class of overdetermined systems of PDEs in \cite{KL$_1$,KL$_2$,KL$_4$}.

We briefly recall the results. If $\E=\{F_1=0,\dots,F_r=0\}$ is a system of PDEs ($r=\dim H^{*,1}(\E)$ is an invariant of the system expressed via Spencer $\delta$-cohomology) on $m$ dependent functions $u^1,\dots,u^m$ of $n$ independent variables $x^1,\dots,x^n$, then we call such system a generalized complete intersection if $m\le r<n+m$ and the characteristic variety and characteristic sheaf over it satisfy certain genericity (transversality) conditions (look for details in loc.sit.). The strict inequality means that the system is overdetermined, but not too overdetermined.

For the coordinate situation, which we treat (see \cite{KL$_4$} for a more general context), the multi-bracket of (non-linear) differential operators is defined as follows.
Let $\ell(G)=(\ell_{1}(G),\dots,\ell_{m}(G))$ be the linearization of an operator $G$ written in components. Then
 $$
\hspace{-0.1in}\{F_{1},\dots,F_{m+1}\}=\dfrac{1}{m!}\hspace{-0.1in}%
\sum_{\alpha\in\mathbf{S}_{m},\beta\in\mathbf{S}_{m+1}}\hspace{-0.1in}\left(
-1\right)  ^{\alpha}\left(  -1\right)  ^{\beta}\ell_{\alpha(1)}(F_{\beta
\left(  1\right)  })\circ\ldots\circ\ell_{\alpha(m)}(F_{\beta\left(  m\right)
})\left(  F_{\beta(m+1)}\right).
 $$
For $m=1$ (scalar differential operators) this becomes the Jacobi bracket $\{\,,\}$ on $\op{diff}(\1,\1)$, which in the linear case become the standard commutator.

With the condition of generalized complete intersection a system $\E$ is compatible (formally integrable) iff the multi-brackets vanish due to the system, meaning
 \begin{equation}\label{multbr}
\{F_{i_1},\dots,F_{i_{m+1}}\}=0\,\op{mod}\mathcal{J}_{k_{i_1}+\dots+k_{i_{m+1}}-1}(F_1,\dots,F_r)
 \end{equation}
for any $(m+1)$-tuple $1\le i_1<\dots<i_{m+1}\le r$, where $k_i=\op{ord}(F_i)$ are orders of the operators and $\mathcal{J}_s(F)=\langle\hat\square_i\circ F_i\,\bigl|\,\op{ord}\square_i+k_i \le s,1\le i\le r\rangle\subset\op{diff}_s(m\cdot\1,\1)$ is the submodule generated by $F_1,\dots,F_r$ and their total derivatives up to order $s$.

Another important feature is the description of the solutions space. In particular, functional dimension and rank of the space $\op{Sol}(\E)$ can be calculated explicitly \cite{KL$_3$}, which helps to control the reduction process, which we shall now discuss.

\subsection{\hpss Symmetries and Lagrange-Charpit method.}\label{S11}

 \abz
Let $F=(F_1,\dots,F_m)^T$ be a determined differential operator (this is not only condition on the number of functions, but also the claim that determinant of the symbolic matrix does not vanish identically). Another operator $G=(G_1,\dots,G_m)^T$ is a symmetry of $F$ if
 \begin{equation}\label{sym}
\{F,G\}=\l\circ F
 \end{equation}
for some matrix differential operator $\l$ (or its linearization in non-linear case)\footnote{Note that $\op{ord}(\l)\le\op{ord}(G)$ for $m>1$ and it is $<\op{ord}(G)$ for $m=1$.}. This condition means that the overdetermined system $F[u]=0,u_t=G[u]$ is compatible ($t$ is an additional independent variable), which can be interpreted as a vector field $\re_G$ on the solution space $\op{Sol}(F)$ \cite{KLV} (flow of the evolutionary equation).

With certain non-degeneracy conditions for the symbols of $F,G$ equation (\ref{sym}) implies that the joint system $F=0,G=0$ is compatible, which can be interpreted as existence of auto-model solutions for $G$ (fixed points of the flow $\re_G$)\footnote{At this point I am grateful to S.Igonin and A.Verbovetsky for an enlightening discussion of the results of \cite{KL$_4$} and the symmetry condition.}. Indeed, if $m=1$ this was noted already in \cite{KL$_1$}: the compatibility condition, specification of (\ref{multbr}), given by
 \begin{equation}\label{comp}
\{F,G\}=\l\circ F+\mu\circ G\qquad\qquad \bigl(\,\op{ord}(\l)<\op{ord}(G),\,\op{ord}(\mu)<\op{ord}(F)\,\bigr)
 \end{equation}
is the direct consequence of (\ref{sym}). For $m>1$ Theorem 1 of \cite{KL$_4$} is not directly applicable, so either a version of it for non-commutative rings should be exploited or the reduction via D-modules to the scalar system should be applied. This will be performed elsewhere.

Now overdetermination of the equation $F=0$ with another PDE $G=0$ (differential constraint) is the subject of most reduction methods (since $m=1$ in most of them, we restrict our attention mainly to this case now). However without transversality condition (complete intersection of \cite{KL$_1$}) the compatibility condition (\ref{comp}) will not work. And in fact the reduction may not work either, see \S\ref{S01}. Thus it is very important to control this general position property (which ad hoc is present in all reduction methods, though this is not specified explicitly).

Classical Lagrange-Charpit method designed to integrate scalar 1st order PDEs consists of overdetermining them in a special way. If $F_1=0$ is a given scalar PDE, the auxiliary PDEs are chosen $F_2=0,\dots,F_n=0$, so that the resulting system is of Frobenius type (i.e. $\{u_{x_i}=F_i(x,u)\},1\le i\le n$) and compatible. Two points at this method are to be underlined.

First is that the condition formulated by them is given in terms of vanishing of Mayer (or Lagrange) brackets due to the system, so that the criterion of \cite{KL$_4$} is the higher vector-version of the classical theorem. The second is that the symbol of the joint system is of general position and the system is not too overdetermined.

With the general compatibility criterion of \cite{KL$_4$} we can skip the Frobenius condition, but we should keep the condition for symbol and that the system is not too overdetermined. The latter means that we try to put a minimal number of differential constraints. Of course, this limits the general idea, but then an efficient formula for checking compatibility exists\footnote{This is not necessary: Provided that the form of symbolic syzygy (type of singularity of the characteristics) is fixed, the compatibility formulas can be specified. This is the essence of the method (not theorem) of \cite{KL$_4$}.}.

This method was called the generalized Lagrange-Charpit method in \cite{KL$_2$} (it works for vector systems as well, with multi-brackets used instead of the Jacobi bracket). The equation $G$ is called then an auxiliary integral (the terms differential constraint and side condition can be applied, but let us stress importance of non-degeneracy condition).

Restricting for simplicity to scalar PDEs (this is what basically considered in all methods), we can state:

 \begin{quote}
{\it Most reduction methods are partial cases of the generalized Lagrange-Charpit method.\/}
 \end{quote}

This is not true for all methods, since certain ansatzes lead to more overdetermined systems and can be interpreted as solvable differential constraints, not compatible differential constraints. But nevertheless this makes a convenient specification of the general Olver's claim from \S\ref{S02}.

\subsection{\hpss Similarity and other reductions.}\label{S12}

 \abz
Here we specify how to get major reduction methods from the one considered in the previous section. Some of the relations described below are known (see \cite{Ol} and the references therein), but our form is most general and is based on the rigorous compatibility results. We restrict to the case $n=2$ (independent variables $x,t$), $m=1$ (dependent variable $u$) and one single PDE (most reduction methods are discussed precisely for this particular case), though versions for the general case are straightforward.

\smallskip

{\bf 1.} {\it Symmetries and generalized symmetries.\/} The point symmetries are quasi-linear auxiliary integrals of the form:
 \begin{equation}\label{quasil}
G=a(t,x,u)u_t+b(t,x,u)u_t+c(t,x,u).
 \end{equation}
The only difference between classical symmetries introduced by S.Lie \cite{Lie,LE} and generalized symmetries of Bluman and Cole \cite{BC}, see also \cite{Ov,Me}, is that they satisfy different relations (\ref{sym}) and (\ref{comp}) respectively, of which the latter is more general. Thus for establishing (partially) invariant solutions one can treat them without any distinction.

Contact symmetries were, also introduced by S.Lie and they are general 1st order constraints $G=G(t,x,u,u_t,u_x)$, which satisfy (\ref{sym}). Though never treated generalized contact symmetries should be the general first order $G$ satisfying (\ref{comp}). They are more plentiful than generalized point symmetries.

\smallskip

{\bf 2.} {\it Higher symmetries.\/} These are just differential operators $G$ of higher order satisfying (\ref{sym}). Upgrading it to (\ref{comp}) we get higher generalized symmetries or auxiliary integrals (as well as compatible differential constraints and side conditions). The only restriction which should not be forgotten in all these overdeterminations is the non-degeneracy condition (i.e. $F=G=0$ is a complete intersection).

Higher order auxiliary integrals exist in abundance, but to find from their continuum certain concrete families one should specify the form of $G$ due to original PDE (see examples in \cite{KL$_1$}).

\smallskip

{\bf 3.} {\it Direct method and beyond.\/}
Direct reduction of Clarkson and Kruskal \cite{CK} is an ansatz for the form (\ref{quasil}), so it gives the same exact solutions as explained in 1. This method was generalized by Galaktionov to the following ansatz: $u=U(t,x,w_1(\zeta(t,x)),\dots,w_k(\zeta(t,x))$ (he took $k=2$) and the substitution is aimed to reduce the PDE to an ODE of higher order for the variable $z=\zeta(t,x)$.

If $V$ is the vector field with first integral $\zeta(t,x)$, then Olver \cite{Ol} demonstrates that this is equivalent to a differential constraint of the form $G=V^k(u)-\Phi(t,x,u,V(u),\dots,V^{k-1}(u))=0$,
which is a special form (degenerate symbol, characteristic variety is $k$-multiple point in $\C P^1$) of quasi-linear PDE.

Provided that differential operator $F$ of order $m$ has $\op{Char}^\C(F)$ not containing this point and that $\{F,G\}=0\,\op{mod}\mathcal{J}_{k+m-1}(F,G)$, the PDE $F=0$ is reduced to an ODE of order $km$ \cite{KL$_1$}.

 \begin{rk}
If an ansatz reduces PDE $F=0$ of order $k$ to an ODE of order $N$ or a system with solutions space of dimension $N$, and number $N$ is not divisible by $k$, then this reduction cannot be obtained by a single compatible differential constraint. Indeed transversality is necessary for reduction and then compatibility implies $\dim\op{Sol}(F,G)=\op{ord}F\cdot\op{ord}G$. Thus this case corresponds to several (different orders) jointly compatible differential constraints or to solvable differential constraints (this can be a single scalar $G$).
 \end{rk}

Generalization of this method to arbitrary $n$ and $m$ is this: One specifies $G$ to have degenerate form and characteristic variety. Darboux integrability then becomes a partial case.

\smallskip

{\bf 4.} {\it Separation ansatzes.\/}
Additive separation ansatz $u=f(t)+g(x)$ (we continue the case $n=2$, $m=1$ for simplicity) correspond to $G=u_{xy}$, multiplicative one $u=f(t)g(x)$ corresponds to $G=uu_{tx}-u_tu_x$.
Functional separation $u=U(f(t)+g(x))$ generalizes these both and corresponds to
$G=u_tu_x^2u_{ttx}-u_t^2u_xu_{txx}+u_{tx}(u_t^2u_{xx}-u_x^2u_{tt})$.

Functional separation of variables $u=\sum_{i=1}^kf_i(t)g_i(x)$ corresponds to a complicated non-linear higher order constraints. Note that for $k>1$ encoding the ansatz cannot be chosen as a single PDE. For instance, when $k=2$ we can have either a mixed ansatz-constraint (i.e. an auxiliary integral $G$ containing an unknown function)
 $$
(u_{tx}u-u_tu_x)p''(t)+ (u_xu_{tt}-u_{ttx}u)p'(t)+ (u_{ttx}u_t -u_{tt}u_{tx})p(t)=0
 $$
or 2 compatible PDEs of order 5, expressing non-trivial solvability of the above equation for $p(t)$.

\smallskip

{\bf 5.} {\it $N$-soliton and finite-gap solutions of soliton equations.\/}
These are always special type constraints. Consider for instance KdV equation $u_t=u_{xxx}+6uu_x$. $N$-soliton solution is obtained from an ansatz\footnote{The parameters should be in fact more specific, otherwise this family yields more solutions.} $u=c\,w_{xx}$, where $e^w=\sum c_ie^{\alpha_i t+\beta_i x}$. Functions $w$ of this type can be described as solution to a PDE $F[w]=e^{-w}L[e^w]=0$, where $L$ is a linear differential operator with constant coefficients. $F$ is a non-linear operator and there exist other non-linear operators of high order $\Phi,\Psi$ such that $\hat\Psi\circ F[w]=\hat\Phi\circ c\,\D_x^2$. Then the (solvable) differential constraint is $\Phi[u]=0$.

The space of $N$-gap solutions (with specified place for gaps) is a finite-dimensional (symplectic) manifold invariant under KdV dynamics. The standard description is however in non-differential terms (spectral properties of the Hill equation), but this can be re-written via differential constraint with $\tau$-functions as coefficients. Since the latter are in turn solutions to differential equations, the auxiliary integrals can be specified in an elementary, but very non-linear form.

Notice also that Painlev\'e test is a special kind of ansatz and so can be seen as a constraint.

\smallskip

{\bf 6.} {\it Evolutionary equations and inertial manifolds.\/}
The generalized Lagrange-Charpit method is much simpler for evolutionary PDEs. Indeed if
$F=u_t-\Phi(t,x,u,u_x,u_{xx},\dots)$, so that no $t$-derivative enters into the last term ($x$ and $u$ can be multidimensional), then an auxiliary integral, with no loose of generality, can be represented in the form $G=\Psi(t,x,u,u_x,u_{xx},\dots)$. Indeed all derivatives of $u_t$ can be expressed from the equation $F=0$.

Seeming difficulty that the characteristic variety of the system $\E=\{F=G=0\}$ is degenerate for 1-dimensional space variable $x$ and $\op{ord}\Phi>1$, can be resolved by introducing higher weight for the variable $t$ \cite{KL$_1$}. If we assume moreover that the auxiliary integral does not contain $t$ explicitly, i.e. $G=G(x,u,u_x,u_{xx},\dots)$, then $G=0$ gives an equation of the stratified submanifold $\mathcal{S}$ (if $n-1=\dim\{(x)\}>1$ it is an infinite-dimensional manifold, while in the case $n=2$ we get an ODE), which is invariant under evolutionary flow, see \cite{KL,LL} for examples.

Such (finite-dimensional) invariant submanifolds $\Sigma=\{G=0\}$ are called inertial if in addition they attract evolutionary dynamics $u_t=\Phi[u]$. In this case $\Sigma$ contains an attractor, which can be investigated by the method of dynamical systems. When exist, inertial submanifolds can be found via compatible differential constraints.

\subsection{\hpss Examples.}\label{S13}

 \abz
Since the method is the most general, we cannot classify all solutions obtained via it (this would be just the description of $\op{Sol}(\E)$), but we show some relations with what other did and something new as well (see also \cite{KL,LL} for other examples).

We will restrict to the case $n=2,m=1$ and one PDE $F=0$. Assuming complete intersection condition, its compatibility with the differential constraint $G=0$ is given by $[F,G]_\E=0$, where the latter bracket is the reduction of the Jacobi bracket $\{F,G\}$ by the subspace $\mathcal{J}_{\op{ord}F+\op{ord}G-1}(F,G)$; this reduction is called the Mayer bracket \cite{KL$_1$}.

{\bf 1.} For $F=\Delta u-1$ (where $\Delta=D_t^2+D_x^2$ is the Laplace operator on $\R^2$) the constraint $G=(\sqrt{u})_t^2+(\sqrt{u})_x^2-\alpha$ is an auxiliary integral iff $\alpha=\frac12$ or $\frac14$. This follows from the formula $[F,G]_\E=\frac{(1-2\alpha)(1-4\alpha)}{2u}$.
 Note that $G$ is non-linear, while $F$ is linear non-homogeneous (cf. \S\ref{S21}).

{\bf 2.} For $F=c(t)u_{tt}-u_{xx}$ we let $G=u_t+a(t)u$. Since $u\cdot[F,G]_\E=c(t)a''(t)+c'(t)a'(t) -2c(t)a(t)a'(t) -a(t)^2c'(t)$, the compatible constraint of the given form corresponds to the function $a(t)=c_0+c_1\int\frac{dt}{c(t)}$. It is not a point symmetry, so it is a non-classical symmetry.

{\bf 3.} Let us re-expose an example from \cite{OR}. For the heat equation $F=u_t-u_{xx}$ they take $G=tu_t-xu_x+3x^3$. This is not a compatible differential constraint, because $H=\frac13[F,G]_\E=u_t-6x$, which is not identically zero. But adding this PDE makes the system compatible. Indeed, $G=H=0$ is a compatible Frobenius system with the solution $u=x^3+6tx+c\in\op{Sol}(F)$.

{\bf 4.} We continue with the heat equation $F=u_t-u_{xx}$. The constraint of the multiplicative ansatz $G=u_{tx}u-u_tu_x$ does not satisfy the complete intersection condition. The compatibility condition is simpler: $H=u\mathcal{D}_t(F)+\mathcal{D}_x(G)-u_tF=u_{tt}u-u_t^2=0$. Thus again
the system $F=H=0$ is not compatible, but the "projected" system $F=G=H=0$ is of Frobenius type and is compatible.

{\bf 5.} Finally let us consider a more complicated case. Choose
 \begin{equation}\label{KPP}
F=u_t-u_{xx}-f(u)=0.
 \end{equation}
This is a non-linear reaction-diffusion equation (called also Kolmogorov-Petrovsky-Piskunov and non-stationary heat equation). The following is an auxiliary integral:
 \begin{equation}\label{KPP1}
G=u_xu_{xxx} -u_{xx}^2-f(u)u_{xx} + f'(u)u_x^2=0.
 \end{equation}
This follows readily from the identity $[F,G]_\E=0$, but can be computed directly from invariance.

Now an important thing is that PDE (\ref{KPP}) was a popular target for investigation for exact solutions. Detailed analysis and a review is contained in \cite{CM}. In all the results $f(u)$ has some particular form, most often cubic. In \cite{KL} this equation was investigated with a different method and still some exact solutions for some other particular form of $f(u)$ were obtained (an exception was with first order auxiliary integrals, when $f(u)$ was general, leading to traveling waves).

But the auxiliary integral (\ref{KPP1}) works for all $f(u)$ and moreover it is already an ODE for $u$ in $x$, so it is the reduction (the solution of the evolutionary equation is reduced to two quadratures: first is $G=0$ and the second is given by integrating the vector field corresponding to $F$ on $\op{Sol}(G)$, see \cite{KL,LL}).

 \begin{rk}
This shows that the method of differential constraints (Lagrange-Charpit in this case) is more powerful than the ansatz approach. The reason is that substitution of separation ansatzes leads to a functional equation, which is difficult to solve in general (indeed following Abel we solve functional equation via differential ones).
 \end{rk}

{\bf 6.} In order to explain the result of {\bf5} let us mention that (\ref{KPP}) was investigated for functional separation of variables (as in direct reduction method) $u=U(\vp(t)x+\psi(t))$ in \cite{PZ}. The corresponding functional equation was partially solved with certain type of $f(u)$ specified.

Now the above ansatz is equivalent to the constraint $\bigl(\dfrac{u_t}{u_x}\bigr)_{xx}=0$, which coupled with (\ref{KPP}) yields a forth order PDE (after multiplication by $u_x^3$)
 $$
G_1=u_x^2u_{xxxx}-(3u_xu_{xx}+f(u)w_x)u_{xxx}+ 2u_{xx}^3+ 2f(u)u_{xx}^2 - f'(u)u_x^2u_{xx}+
f''(u)u_x^4=0.
 $$
However it is not compatible with (\ref{KPP}): $[F,G]_\E=2\,G\cdot G_2$, where $G$ is given by (\ref{KPP1}) and $G_2=f(u)u_xu_{xxx}- 2f(u)u_{xx}^2+
 f'(u)u_x^2u_{xx} -f''(u)u_x^4$.

So the system splits. The first branch $F=G=0$ is compatible and was noted above. This choice gives unconstrained $f(u)$.

The other system $F=G_2=0$ is not compatible. Compatibility analysis implies that the system
$F=G_2=0$ is solvable iff $f=f(u)$ satisfies the following highly nonlinear ODE:
 \begin{multline*}
2f'^4(f'')^2f^{(4)} + f'^3f''f^{(3)}( 7(f'')^2 + 12ff^{(4)} )  +
  ff'f^{(3)}\Bigl( 30(f'')^4 + 45ff''(f^{(3)})^2 + 16f(f'')^2f^{(4)} \\
 + 12{f}^2(f^{(4)})^2 \Bigr) +f'^2\Bigl( 3(f'')^5 + 39f(f'')^2(f^{(3)})^2 + 3f(f'')^3f^{(4)} + 18{f}^2(f^{(3)})^2f^{(4)} + 4{f}^2f''(f^{(4)})^2 \Bigr)  \\
 + f\Bigl( 8(f'')^6 + 31f(f'')^3(f^{(3)})^2 + 45{f}^2f''(f^{(3)})^2f^{(4)} -
     6{f}^2(f'')^2(f^{(4)})^2 + {f}^2( -27(f^{(3)})^4 + 2f(f^{(4)})^3 )  \Bigr)=0.
 \end{multline*}

Thus in this case we also get a reduction, but again for a special kind of $f(u)$.

{\bf 7.} Note that the operator in (\ref{KPP}) is $F=\p_t-L$ with $L[u]=u_{xx}+f(u)$ and $G=\nabla\circ L$ with $\nabla = u_xD_x-u_{xx}$. Thus the evolutionary flow of $F$ is given by the second order ODE on the 3-dimensional space $\op{Sol}(G)$ of solutions of ODE $G=0$ in $x$-variables only. Note that $\op{Sol}(L)\subset\op{Sol}(G)$ is the 2-dimensional surface, on which the flow is stationary. On the complement the flow is non-trivial. This geometric reasoning shows that while there are many auxiliary integrals of the form $G=\square\circ L$, not every such $G$ is a compatible constraints, but any is a solvable differential constraint (indeed $\op{Sol}(L)$ will be a part of $\op{Sol}(F,G)$).

Let us also note that the solutions of (\ref{KPP})+(\ref{KPP1}) can be obtained by a Frobenius overdetermination. Indeed it is equivalent to the joint system
 $$
u_{tt}=\frac{u_t^2}{u_x^2}(u_t-f(u)),\quad u_{tx}=\frac{u_t}{u_x}(u_t-f(u)),\quad u_{xx}=u_t-f(u),
 $$
which is compatible ($F=0$ is the last equation).

\subsection{\hpss Algebraic structure of the set of auxiliary integrals.}\label{S14}

 \abz
In many papers devoted to differential constraints (for instance \cite{Ol,Cl}) it is noted that the method is too general, so it is impossible to operate with, contrary to the symmetry group approach. It is indeed quite general due to Theorem \ref{thm1}. But being reasonably restricted it remains quite general/powerful, while will possess a certain algebraic structure\footnote{In too general form, we observe no algebraic structure for compatible differential constraints. The same concerns solvable constraints, even if the form is prescribed.}, which we are going to explain now.

For simplicity we restrict to the case of a single scalar PDE $F=0$, where $F\in\op{diff}(\1,\1)$ is a (non-linear) differential operator over $M$, as this is an often target of investigations (but change of the bundle $\1$ to $\pi$ is possible with some assumptions on the operators).

 \begin{rk}
Non-classical symmetries are special kind of auxiliary integrals in generalized Lagrange-Charpit method, so their totality does not possess algebraic structure as well. Sometimes to stress group approach it is said that generalized symmetry $G$ of $F=0$ can be characterized by the property that the flow of it preserves the solutions $\op{Sol}(F,G)$ (contrary to the larger space $\op{Sol}(F)$ for classical symmetries). But then general auxiliary integrals are the same kind of generalizations of the higher symmetries (we can treat the flow as defined  formally).
 \end{rk}

In what follows we write $\ell_F=\sum F_{p_\z}\D_\z$ the operator of linearization. Here $(x^i,p_\z)$ are the canonical coordinates on the jet-space $J^\infty(M)$ induced by a coordinate chart on the base $M^n$, $\z=(i_1,\dots,i_n)$ is the multi-index and $\D_\z=\D_1^{i_1}\cdots\D_n^{i_n}$ is the higher total derivative operator with $\D_i=\p_{x^i}+\sum p_{\z+1_i}\p_{p_\z}$. Thus the operator of universal linearization transforms non-linear operators $\op{diff}(\1,\1)$ to linear operators
$\op{Diff}(\1,\1)$ if we fix a jet $x_\infty\in J^\infty(\1,\1)$.
More generally it is the map
  $$
\ell:\op{diff}(\1,\1)\to\Cc\op{Diff}(\1,\1)=C^\infty(J^\infty(\1))\otimes_{C^\infty(M)}\mathop{\rm Diff}(\1,\1)
 $$
(the latter are called scalar $\Cc$-differential operators).

Let us denote $\op{Sym}(F)=\{H\in\op{diff}(\1,\1):\{F,H\}=\ell_\theta F\text{ for some }\theta\in\op{diff}(\1,\1)\}$ the space of symmetries of $F$, which is a Lie algebra w.r.t. the Jacobi bracket. The defining relation determines the family of vector subspaces $\op{Sym}_\theta(F)=\{H\in\op{diff}(\1,\1):\{F,H\}=\ell_\theta F\}$. Note that $\op{Sym}_0(F)$ is the Lie subalgebra of $\op{Sym}(F)$ and $\op{Sym}_\theta(F)$ is the module over it (considered as a vector space!).

Notice that the following space is a vector space, but in general it is not a Lie algebra:
 $$
\op{sym}^*(F)=\{\theta\}=\op{Sym}(F)/\op{Sym}_0(F).
 $$

Next for $\mu\in\op{diff}(\1,\1)$ with $\op{ord}(\mu)<\op{ord}(F)$ denote
 $$
\op{Aux}_\mu(F)=\{G\in\op{diff}(\1,\1):\{F,G\}=\ell_\l F+\ell_\mu G\text{ for some }\l\in\op{diff}(\1,\1)\}
 $$
and let $\op{Aux}(F)=\cup_{\mu:\op{ord}(\mu)<\op{ord}(F)}\op{Aux}_\mu(F)$ be the set of all auxiliary integrals. Notice that we omitted here the condition of general position for the symbols of $F$ and $G$ (complete intersection). With this condition the joint system $\{F=G=0\}$ is compatible and $G$ is a compatible differential constraint.

$\op{Aux}_\mu(F)$ is a vector space.
Changing $G\mapsto G+\ell_\eta F$ we change coefficient $\l$ in $\{F,G\}$. The equivalence classes form the set $\op{aux}^*(F)=\{\mu\}$, which in general is just a topological space.

It's easy to see that $\op{Aux}(F)$ is neither a vector space nor a module over $\op{Sym}(F)$ (not even over $\op{Sym}_0(F)$ as vector space). Indeed, if $G\in\op{Aux}(F)$, $\{F,G\}=\ell_\l F+\ell_\mu G$, and $H\in\op{Sym}_\theta(F)$, then $\{F,G+H\}=\ell_{\l+\theta} F+\ell_\mu G$. Only $\op{Aux}_0(F)=\op{Sym}(F)$ is the module.

However it has the structure of algebroid type. To describe it we need to introduce some notations. Denote $\op{ad}_H=\{H,\cdot\}$. Let $\re_G$ be the evolutionary differentiation, which is the dual operator to universal linearization \cite{KLV}: $\re_GF=\ell_FG$. Note that $\op{ad}_H=\ell_H-\re_H$.

The Hessian operator $\op{Hess}:\op{diff}(\1,\1)\times\op{diff}(\1,\1)\to\Cc\op{Diff}(\1,\1)$ is defined by the formula $\op{Hess}_FG=[\re_G,\ell_F]$. In canonical coordinates we have:
 $$
\op{Hess}_FG(H)=\op{Hess}_F(G,H)=\sum F_{p_\z p_\tau}\D_\z(G)\D_\tau(H).
 $$
For linear differential operators $F$ the operator $\op{Hess}_F$ vanishes. These operators satisfy the following properties (direct check):

 \begin{prop}
1. Symmetry: $\op{Hess}_F(G,H)=\op{Hess}_F(H,G)$.\\
2. Compensated Leibniz rule: $\{F,\ell_GH\}=\ell_{\{F,G\}}H+\ell_G\{F,H\}-\op{Hess}_F(G,H)$.\\
3. Linearization anomaly: $[\ell_F,\ell_G]=\ell_{\{F,G\}}+(\op{Hess}_GF-\op{Hess}_FG)$.
 \end{prop}

Now we can use this to calculate:
 \begin{multline*}
 \{F,\{G,H\}\}= \{\{F,G\},H\}\}+\{G,\{F,H\}\}
= -\{H,\ell_\lambda F+\ell_\mu G\}\}+\{G,\ell_\theta F\}\\
= (\ell_{\{\lambda,H\}}+[\ell_\lambda,\ell_\theta]-\ell_{\{\theta,G\}}+\op{Hess}_H\lambda-
\op{Hess}_G\theta)F+\ell_\mu\{G,H\}\\
 +(\ell_{\{\mu,H\}}-\ell_\theta\ell_\mu+\op{Hess}_H\mu)G
\ =\ -(\op{ad}_H+\ell_\theta)\ell_\mu G\ \op{mod}F.
 \end{multline*}

From the last two lines we conclude two actions:

 \begin{prop}
1. If $H\in\op{Sym}_\theta(F)$, $(\ell_{\op{ad}_H}+\ell_\theta-\op{Hess}_H)\mu=0$, then $\op{ad}_H:\op{Aux}_\mu(F)\to\op{Aux}_\mu(F)$.\\
2. If $H\in\op{Sym}_\theta(F)$, i.e $(\op{ad}_H+\ell_\theta)F=0$, and
$(\op{ad}_H+\ell_\theta)\ell_\mu=0$, then $\op{ad}_H:\op{Aux}_\mu(F)\to\op{Sym}(F)$.
 \end{prop}

Thus we get partial algebraic structure on the set $\op{Aux}(F)$, so that the fibers of the vector bundle $\op{Sym}(F)\to\op{sym}^*(F)$ act on the fibers of the vector bundle $\op{Aux}(F)\to\op{aux}^*(F)$.

\section{\hps Other related results}\label{S2}

 \abz
The compatibility criterion via multi-brackets discussed in the previous section can be generalized to the case, when we still have $r<m+n$ equations on $m$ dependent functions, but allow degenerations keeping regularity. Let us begin for simplicity with the case $m=1,r=2$. For complete intersection we required that the characteristic varieties $\op{Char}^\C(F)$ and $\op{Char}^\C(G)$ are transversal in each fiber $\C P^{n-1}$ (\cite{KL$_1$}). We can relax this by letting them have common components. Assume the symbols $\z(F),\z(G)$ have a (greatest) common divisor $q$. Then we can reduce the bracket by an operator $Q$ with $\z(Q)=q$.

In fact, let $k=\op{ord}(F)$, $l=\op{ord}(G)$, $t=\op{ord}(Q)$ and $S,T$ be differential operators with symbols $\z(F)/q$, $\z(G)/q$ (regularity assumption means the division is made over all points of the system $\E=\{F=0,G=0\}$ and that loci of polynomials $\z(F)/q,\z(G)/q$ intersect transversally). Then $\E$ is formally integrable iff $\ell_TF-\ell_SG=0\,\op{mod}\mathcal{J}_{k+l-q-1}(F,G)$. The operator on the left is the reduced analog of the bracket.

Generalization to the case $2<r\le n$ is straightforward. Multi-bracket analog $m>1$ is more involved. It is also possible to have explicit criteria when $r\ge m+n$, but this will be treated elsewhere.
In this section we consider some questions related to compatibility issues of establishing exact solutions.

\subsection{\hpss ODEs.}\label{S2a}
 \abz

In \cite{Cl}\S6.3 it was noticed that contrary to classical symmetry reduction the nonclassical and direct method are difficult to use for ODEs. In fact, there's no difficulties if we view them as differential constraints. However in this case the number of independent variables cannot be reduced (it is already 1) and the order should be reduced instead (cf. Remark \ref{rk2}).

Consider an example, which has a general nature. Let
 \begin{equation}\label{ode1}
u''=F(x,u,u')
 \end{equation}
be an ODE. Seek for the compatible differential constraint in the form $u'=\vp(x,u)$. In this case the second equation is just the intermediate integral of the first equation and condition of this is:
 \begin{equation}\label{ode2}
\vp_x=-\vp_u\vp+F(x,u,\vp).
 \end{equation}
Thus reductions of ODE (\ref{ode1}) are bijective to partial solutions of PDE (\ref{ode2}). This is an analog of Hamilton-Jacobi approach.

It is not clear how to find a solution to PDE, but whenever it is done we proceed to simplify (\ref{ode1}). For instance all classical reductions can be easily recognized. Let's mention the basic:

1) If $F$ is independent of $x$, then $x$ does not enter (\ref{ode2}) and we seek for solution $\vp=\vp(u)$. Then (\ref{ode2}) is the ODE-reduction.

2) If $F$ is linear in $u$, then from (\ref{ode2}) we read that it worth searching $\vp$ to be homogenous in $u$ of degree 1.

3) If $F(e^{\alpha t}x,e^{\beta t}u,e^{(\beta-\alpha)t}\vp)=e^{(\beta-2\alpha)t}F(x,u,\vp)$, one can take $\vp(e^{\alpha t}x,e^{\beta t}u)=e^{(\beta-\alpha)t}\vp(x,u)$.

In all these cases one uses a symmetry of (\ref{ode1}) to be the symmetry for (\ref{ode2}), which in the considered cases are respectively: $\p_x$, $u\p_u$ and $\alpha x\p_x+\beta u\p_u$.

On the other hand what we do is just to solve (\ref{ode2}) with the help of symmetries, which we can write as auxiliary equations:
 $$
1)\ \vp_x=0,\qquad  2)\ u\vp_u-\vp=0,\qquad 3)\ \alpha\vp_x+\beta\vp_u+(\alpha-\beta)\vp=0.
 $$
But in general to solve (\ref{ode2}) we can use more general Lagrange-Charpit method with some auxiliary equations $\vp_u=G(x,u,\vp,\vp_x)$, so that the procedure of solving (\ref{ode1}) becomes repetitive.

\subsection{\hpss Linear systems.}\label{S21}

 \abz
Linear PDEs, as the simplest equations, produce sometimes confusion in problems related to compatibility. Here we discuss some.

{\bf 1.} Olver claims in Theorem 13 of \cite{Ol} that {\it a compatible differential constraint for a second order linear equation must be linear\/}. This is wrong.

A simple counter-example can be obtained for the most popular equation of the sort: heat equation $F=u_t-u_{xx}$. Then $G=(u_t^2+u_x^2)\op{cosh}(2x)-2u_tu_x\op{sinh}(2x)+2e^{2t}(u_t-u)$ is an auxiliary integral, i.e. the system $\E=\{F=G=0\}$ is compatible.

The order of $G$ is 1 and $\op{Sol}(\E)=\{c_1e^{t-x}+c_2e^{t+x}+c_1^2+c_2^2\}$. The point is that this is a curved surface in the linear space $\op{Sol}(F)$, in fact a quadric. By this reason the function $G$ is quadratic.

Similarly if we take a 4-dimensional curved surface in $\op{Sol}(F)$, we get an auxiliary non-linear integral $G$, i.e. a compatible constraint contradicting Olver's claim (he was concerned with second order nonlinear constraints). For instance the following auxiliary integral works:
 \begin{multline*}
G=x^2(2u_t^2 + 2u_{tt}^2-u_{tt}+u_t) +2u_{tt}( 1- t+2xu_{tx} -2xu_x)\\  +
  2u_t( t - 2x^2u_{tt} - 2xu_{tx} + 2xu_x) +2u_{tx}(u_{tx}-x - 2\,{u_x} )+ 2(xu_x + {u_x}^2-u). \end{multline*}

It is quite clear that the construction is most general and plenty of nonlinear constraints for linear equations exist\footnote{Constant coefficient linear PDEs are often target of investigation for particular solutions via reduction methods to excitement of applied mathematicians, who find {\it all solutions\/} by the Fourier transform.\label{ftn12}}.
For instance $F=u_{tt}\pm u_{xx}$ possesses the compatible differential constraint $G=\op{det}\op{Hess}(u)+ae^{bu}$ with Monge-Ampere type non-linearity (see \cite{KL$_1$}).

{\bf 2.} For linear PDEs with constant coefficients an auxiliary integral can be chosen of the same form: A system of PDEs with constant coefficients is compatible, provided it is of generalized complete intersection type, as follows from the compatibility criterion of Section \ref{S1}. Without this requirement the claim is wrong.

For instance the system $\{u_{tx}=u,u_t=u\}$ has only trivial solution ($u=0$ is always a solution of a linear homogeneous system). Also for linear non-homogeneous system with constant coefficients the above claim would be wrong: The system $\{u_t=u+1,u_x=1\}$ has no solutions.

Since any two linear operators $L_1,L_2$ with constant coefficients commute $\{L_1,L_2\}=0$, one is a symmetry for the other. But without complete intersection condition the PDEs $L_1[u]=0$, $L_2[u]=0$ might be incompatible. This shows importance of our generic position assumption.

{\bf 3.} Investigating whether nonclassical method gives more symmetric solutions than the classical one, when applied to linear PDEs, the authors of \cite{AGB} make the following statement:

 \begin{quote}
{\it Every solution of a linear PDE is invariant under some classical symmetry.\/}
 \end{quote}

Then they prove it under the following assumptions: a) $n=2$ (two independent); b) $m=1$ (one dependent); c) the PDE has constant coefficients.

It is a bit disappointing that the proof contains determining relations for symmetries, calculations etc. As the authors expect that their theorem can be extended to bigger number $n$, we confirm this and make the following generalization: We will consider systems of $r$ equations on $m$ unknowns and will not require that coefficients of the PDEs be constant.

 \begin{prop}
Let $L\in\op{Diff}(m\cdot\1,r\cdot\1)$ be a linear operator. Then every solution to the system $L[u]=0$ is invariant under some symmetry provided the following condition on the symmetry algebra holds:
 $$
\op{Sym}(L)\supsetneqq \op{Sym}_\flat(L)\simeq\op{Sol}(L).
 $$
 \end{prop}

Here we denote by $\op{Sym}_\flat(L)\subset\op{Sym}_0(L)$ the subalgebra of shifts by solutions. The symmetries have the form $g\p_u:=\sum g^j(x)\p_{u^j}$ for $g=(g^1,\dots,g^m)\in\op{Sol}(L)$. Thus such a field has the generating function lifted from the base \cite{KLV} and its prolongation is $\sum D_\z(g^j)\p_{u^j_\z}$. The inequality above holds true for operators $L$ with constant coefficients ($\p_{x^1}$ is a symmetry, so that we can allow dependence on $x^2,\dots,x^n$) and thus the claim contains Theorem 1 of \cite{AGB}.

\medskip

 \begin{proof}
Let $u$ be a solution and $\xi$ a symmetry from $\op{Sym}(L)\setminus\op{Sym}_\flat(L)$. Denote $g=\xi u-u$. Then $g\in\op{Sol}(L)$ and the symmetry $\xi-g\p_u$ leaves $u$ invariant.
 \end{proof}

The large size of symmetry algebra does not help in solving linear PDE because the piece $\op{Sym}_\flat(L)$ is generally untractable. One should rather consider the finite piece $\op{Sym}(L)/\op{Sym}_\flat(L)$ (it is given by linear overdetermined system and often is finite dimensional).

Thus even though the conclusion of \cite{AGB} is that no generalized symmetry gives new solutions compared to classical symmetries, this claim does not bear a constructive obstruction for considering generalized symmetries (compatible constraints). For me though footnote \ref{ftn12} does, when PDEs with constant coefficients are considered.

\subsection{\hpss Conclusion.}\label{S22}

 \abz
In this paper we demonstrated how one can establish new exact solutions of non-linear PDEs. The scheme we describe is general enough to include most practiced methods. The method is indeed classical, but we furnished it with an efficient underlying machinery of investigation for compatibility and precise general position conditions to hold.

The generalized Lagrange-Charpit method we discuss was essentially developed in \cite{KL$_2$} and here we make a modification, allowing non generalized complete intersections, but with a fixed type of the syzygy module (then a criterion different from that of \cite{KL$_1$} may show up), or changing compatibility claim to solvability, of which vanishing of multi-brackets due to the system \cite{KL$_4$} is the first and most complicated step.

Thus the proposed technique is not too general, but manageable. We have not performed too many calculations in this paper (see \cite{KL$_1$,KL,LL}), but frankly explained how to do it, leaving more results for future publications. In addition, the method allows to satisfy physical pre-requirements for the desired solutions (form of the compatibility constraint), which might be too restrictive with symmetry approaches.


\vspace{-10pt} \hspace{-18pt} {\hbox to 12cm{ \hrulefill }}
\vspace{2pt}

{\footnotesize \hspace{-20pt} Institute of Mathematics and
Statistics, University of Troms\o, Troms\o\ 90-37, Norway.\ \, E-mail: kruglikov\verb"@"math.uit.no.\\
Gratitude to the research stay in Max Planck Institute for Mathematics in the Sciences, Leipzig, in April-May 2007.}

\end{document}